\documentclass[11pt,a4paper]{article}
\usepackage[T1]{fontenc}
\usepackage[utf8]{inputenc}
\usepackage[english]{babel}
\usepackage[margin=27mm]{geometry}
\usepackage{amsmath,amssymb,amsthm,mathtools}
\usepackage{dsfont,graphicx}
\usepackage{tikz-cd}
\newcommand{\mathbbm}[1]{\mathds{#1}}
\usepackage{tocloft}
\usepackage[colorlinks=true,linkcolor=blue,citecolor=blue,urlcolor=blue]{hyperref}
\newtheorem{thm}{Theorem}[section]
\newtheorem{prop}[thm]{Proposition}
\newtheorem{lemma}[thm]{Lemma}
\newtheorem{example}[thm]{Example}
\newtheorem{construction}[thm]{Construction}
\newtheorem{remark}[thm]{Remark}
\theoremstyle{definition}
\newtheorem{defn}[thm]{Definition}

\newcommand{\Sp}{\mathrm{Sp}}
\newcommand{\SpG}{\mathrm{Sp}^{G}}

\newcommand{\Sph}{\mathbb S}
\newcommand{\Fp}{\mathbb F_p}

\newcommand{\Pp}{\mathcal P}
\newcommand{\Tt}{\mathcal T}
\newcommand{\EP}{E\mathcal P}
\newcommand{\tEP}{\widetilde E\mathcal P}
\newcommand{\Res}{\operatorname{Res}}
\newcommand{\Ind}{\operatorname{Ind}}

\newcommand{\invlim}{\mathop{\underleftarrow{\mathrm{lim}}}}

\newcommand{\NormG}{N_e^G}
\newcommand{\op}{\mathrm{op}}
\newcommand{\coAlg}{\operatorname{coCAlg}}
\newcommand{\NormSys}{\operatorname{NormSys}}
\newcommand{\CAlg}{\operatorname{CAlg}}
\newcommand{\Map}{\operatorname{Map}}
\newcommand{\Fun}{\operatorname{Fun}}
\newcommand{\Span}{\operatorname{Span}}

\newcommand{\Corr}{\operatorname{Corr}_p}
\newcommand{\End}{\operatorname{End}}
\newcommand{\id}{\mathrm{id}}
\newcommand{\infl}{\operatorname{inf}}
\newcommand{\Vect}{\mathrm{Vect}^{\mathrm{fd}}_{\Fp}}
\newcommand{\Kpart}{K^{\mathrm{part}}}

\makeatletter
\DeclareRobustCommand{\lneqslant}{%
  \mathrel{\mathpalette\lneqslant@{}}}
\newcommand{\lneqslant@}[2]{%
  \begingroup
  \setbox0=\hbox{$\m@th#1\leqslant$}%
  \ooalign{\copy0\cr
    \hidewidth\raisebox{-.10\wd0}[0pt][0pt]{%
      \rotatebox[origin=c]{55}{\rule{.42\wd0}{.05\wd0}}}%
    \hidewidth\cr}%
  \endgroup}
\makeatother

\title{The generalized Segal conjecture and the coalgebraic Frobenius}
\author{Thomas Nikolaus}
\date{\today}
\begin{document}
\maketitle

\begin{abstract}
We prove the generalized Segal conjecture, extending work of Bergsaker
and Rognes and confirming a conjecture by the author:
for every finite group $G$ and every bounded-below spectrum $X$,
the Hill--Hopkins--Ravenel norm $N_e^G X$ becomes Borel-complete
after completion at the augmentation ideal of the Burnside ring.
Equivalently, for every non-trivial finite $p$-group $G$, the proper
Tate diagonal $X\to(X^{\otimes G})^{\varphi G}$ is a $p$-completion.
The proof reduces to the case $X=\Fp$, established by Bergsaker
and Rognes, and removes the finiteness hypotheses from their
theorem by abstract manipulations. 

As an application of this result, we refine the coalgebraic
Frobenius to an oriented circle action (i.e. an action of the monoidal category $B\mathbb N$) on the category of bounded-below, $p$-complete 
$\mathbb{E}_\infty$-coalgebras, confirming another conjecture by the author.
The construction is based on Yuan's approach to Frobenius
actions via partial algebraic $K$-theory.
\end{abstract}

\tableofcontents

\section{Introduction}

The classical Segal conjecture, now a theorem of Carlsson~\cite{Car}, is one of the deepest and most mysterious statements in stable homotopy theory. In the last years there has been quite some work on conceptually understanding and (re)proving versions of it. Some of these results, which we will discuss below, have led the author to formulate the generalized Segal conjecture in the ICM proceedings \cite[Conjecture 7.10]{Nik}. In this paper we prove this conjecture using crucial work of Bergsaker--Rognes.

\begin{thm}\label{Tate_diagonal}
For every prime $p$, every non-trivial finite $p$-group $G$ and every bounded-below spectrum $X$, the proper Tate diagonal
\[
 X\xrightarrow{\Delta_G}(X^{\otimes G})^{\varphi G}
\]
exhibits the target as the $p$-completion of $X$. Here $G$ acts on $X^{\otimes G} = X \otimes \ldots \otimes X$ by permuting the factors via left multiplication on $G$ and $(-)^{\varphi G}$ denotes the proper Tate construction.
\end{thm}

We will prove this statement via a (seemingly) more general completion formulation for arbitrary finite groups, as
suggested in \cite[Remark 7.11]{Nik}. We will formulate this using the language of equivariant stable homotopy theory, notably the language of genuine $G$-spectra and the Hill--Hopkins--Ravenel (HHR) norm~\cite[\S2.2.3]{HHR}.\footnote{For earlier multiplicative norm maps, see Greenlees--May~\cite[\S\S4--7]{GM}.} 
Let $A(G)$ be the Burnside ring and 
$I(G)=\ker(A(G)\to\mathbb Z)$ be its augmentation ideal.

\begin{thm}\label{thm:main}
For every finite group $G$ and every bounded-below spectrum $X$, the
$I(G)$-completion of the HHR norm $\NormG X$ is Borel-complete. In particular we have the equivalence of ordinary spectra
\[
(X^{\otimes G})^{hG} \simeq \left ((\NormG X)^G\right)^\wedge_{I(G)} \ .
\]
\end{thm}

Theorems~\ref{Tate_diagonal} and~\ref{thm:main} are equivalent.
In Section~\ref{sec:equivalent}, we explain this equivalence
and the $I(G)$-completion appearing here in detail.

\subsection*{Previously known cases}

For $X=\Sph$, the norm $\NormG\Sph$ is the genuine $G$-sphere.
Thus Theorem~\ref{thm:main} specializes to the classical Segal conjecture
for finite groups, proved by Carlsson~\cite{Car},
building on important earlier work~\cite{Lin,Gun,AGM,MM}.
In degree zero, the classical conjecture identifies the zeroth stable cohomotopy
of $BG$ with the completion of the Burnside ring $A(G)$ at its augmentation ideal.

For $G=C_p$, the proper Tate construction agrees with the ordinary Tate
construction. Lun\o e--Nielsen and Rognes~\cite[Theorem 5.13]{LNR} proved the
norm formulation of the generalized Segal conjecture in this case under the additional assumption that
$H_*(X;\Fp)$ is finite-dimensional in each degree.
The author and Scholze~\cite[Theorem III.1.7]{NS} subsequently removed this
finiteness assumption by reducing to the case of $X = \Fp$ which is covered by \cite{LNR}, establishing the theorem for every
bounded-below spectrum. Further proofs of the case $X=\Fp$ were given in
\cite{HW,Carrick} for $p=2$ and in \cite{AKZ} for odd primes and recently  by the author~\cite{NikSegal} for all primes. Thus the case $G = C_p$ of the two theorems is already known. 
The case $G=C_{p^n}$ follows formally from the case
$G=C_p$ using the Tate orbit lemma, see
\cite[discussion following Conjecture~7.10]{Nik}.

Bergsaker--Rognes~\cite[Theorem 1.1]{BR} proved the norm formulation for all
finite $p$-groups under the same finiteness assumption on mod~$p$ homology.
They also established the $I(G)$-completion formulation for arbitrary finite groups
when $\pi_*(X)$ is finitely generated in each degree \cite[Theorem 1.2]{BR}.
In the present paper we remove these finiteness assumptions.
Our proof builds on their isotropy-separation reduction
\cite[Proposition~2.2]{BR} and uses their $p$-group theorem for $X=\Fp$
as its input. We extend this case to all bounded-below spectra using
exactness, Postnikov towers, and preservation of products on
$\Fp$-vector spaces, following the strategy for $G=C_p$ in
\cite[Theorem~III.1.7]{NS}. The isotropy defect used below is precisely
the obstruction studied by Bergsaker--Rognes.

We hope to also find an independent proof of the case $X = \Fp$ along the lines of \cite{NikSegal} in future work, but for the moment have not succeeded in doing so. 

While preparing this paper, we learned from Vignesh Subramanian
that he had independently obtained a proof of the generalized
Segal conjecture. We also became aware of a note by Wei Yang
\cite{Yang}, dated July 28, 2026, which contains another proof.

\subsection*{Application: The coalgebraic Frobenius}

As an application of these results, we also prove 
\cite[Conjecture 4.5]{Nik} from the author's ICM address, already envisioned as a potential application there (and also by Yuan in unpublished work, see \cite[\S1.6]{Yuan}). To state this result precisely 
recall that an $\mathbb{E}_\infty$-coalgebra in $p$-complete spectra
consists of a $p$-complete spectrum $C$ with comultiplication
$C\to C\widehat{\otimes}_p C$ and counit $C\to\Sph_p^\wedge$,
which are coherently coassociative, counital and cocommutative.
We write $\coAlg_p$ for their $\infty$-category and
$\coAlg_p^{\mathrm{bb}}$ for the full subcategory of coalgebras
whose underlying spectra are bounded below.
For $C\in\coAlg_p^{\mathrm{bb}}$, the \emph{coalgebraic Frobenius}
is the composite
\begin{equation}\label{eq:coalgebra-frobenius}
 \varphi_p:C\xrightarrow{\Delta^{p}}(C^{\otimes p})^{hC_p}
 \xrightarrow{\mathrm{can}}(C^{\otimes p})^{tC_p}
 \xrightarrow{\Delta_{C_p}^{-1}}C.
\end{equation}
Here tensor powers are $p$-completed, the cyclic group acts
by permuting the factors, $\Delta^{p}$ is the iterated comultiplication,
and $\mathrm{can}$ is the canonical map from homotopy fixed points
to the Tate construction. The inverse of the Tate diagonal exists by
the $C_p$-case of Theorem~\ref{Tate_diagonal}, already proved in
\cite[Theorem III.1.7]{NS}. This is the construction of
\cite[Definition 4.3]{Nik}. 


Conjecture 4.5 of \cite{Nik} says that this endomorphism refines to an `oriented circle action' on the category of coalgebras. 
More precisely let $B\mathbb N$ denote the category with one object and
endomorphism monoid $(\mathbb N,+)$, equipped with the symmetric
monoidal structure given by addition. We call it the
\emph{oriented circle}. For coalgebras with $p$-adically finite
underlying spectra, a Frobenius action of this category follows
by duality from Yuan's theorem for algebras
\cite[Theorem 6.8]{Yuan}. Combining his strategy with
Theorem~\ref{thm:main} removes this finiteness assumption:

\begin{thm}\label{thm:coalgebra-frobenius-action}
For every prime $p$, the coalgebraic Frobenius refines to an action
of the monoidal category $B\mathbb N$ on $\coAlg_p^{\mathrm{bb}}$.
More precisely, there is a strong monoidal functor
\[
 B\mathbb N\longrightarrow
 \bigl(\Fun(\coAlg_p^{\mathrm{bb}},\coAlg_p^{\mathrm{bb}}),
 \circ\bigr)
\]
sending the unique object to the identity functor and
$1\in\mathbb N$ to the natural transformation $\varphi_p$.
\end{thm}

As explained in \cite[Remark 4.6]{Nik} and
\cite[Remark 1.14]{Yuan}, the additional content of a $B\mathbb N$-action over a natural endomorphism of each object is
the coherent commutativity of the Frobenius with itself:
the map from $\mathbb N$ to endomorphisms of the identity functor
refines to a map of $\mathbb{E}_2$-monoids. The proof of
Theorem~\ref{thm:coalgebra-frobenius-action} occupies
Part~\ref{part:frobenius}, where we use the generalized Segal conjecture
to adapt Yuan's construction to coalgebras.
The $B\mathbb N$-action constructed here will be used in
forthcoming joint work with Riedel generalizing Yuan's and Mandell's theorems to non-finite settings. 

A $B\mathbb N$-action on all bounded-below $p$-complete
$\mathbb E_\infty$-coalgebras was also announced by Yuan
\cite[Remark~1.5 and \S1.6]{Yuan}. After finishing a first draft of this paper, I learned the details of his unpublished work on the coalgebraic statement, using global coalgebras and a relative form of the generalized Segal conjecture.

\subsection*{Notation}
Throughout, fibers and pullbacks are understood as (homotopy) fibers and pullbacks.
We use the same notation for an abelian group and its associated Eilenberg--MacLane spectrum. Tensor products without further qualifier are always tensor products of spectra (i.e. relative to the sphere $\mathbb S$). 

 We let $G$ be a finite group and write $\SpG$  for the $\infty$-category of
genuine $G$-spectra (see e.g.\ \cite[\S2]{Yuan} for a review). We refer to those as genuine $G$-spectra or simply $G$-spectra. By a spectrum with $G$-action we simply mean an object of $\Sp^{BG}$.
We use the usual notation $Y^H$ and $\Phi^H Y$ for the fixed points and geometric fixed points of a genuine $G$-spectrum $Y$. 
Unless otherwise stated, geometric fixed points are viewed as ordinary spectra. We let $\Pp$ be the family of proper subgroups of $G$ and
$\Tt=\{e\}$ the trivial family. The classifying $G$-anima of these families are denoted $E\Pp$ and $E\Tt = EG$ as usual.

We let $F(-,-)$ denote the inner hom in $G$-spectra. For a pointed $G$-anima $K$, we write
$F(K,Y)$ for $F(\Sigma^\infty K,Y)$ to shorten notation. We write 
$Y^\wedge_{\Pp}=F(\EP_+,Y)$ and $Y^\wedge_{\Tt}=F(EG_+,Y)$
for $\Pp$-completion and Borel completion. For $G\ne e$, the inclusion
$\Tt\subseteq\Pp$ gives natural maps
$
 Y\longrightarrow Y^\wedge_{\Pp}\longrightarrow Y^\wedge_{\Tt}.
$
In this case, the first induces equivalences on proper fixed points and the second on the underlying spectrum. The category of Borel-complete $G$-spectra is equivalent to $\Sp^{BG}$. 

For a prime $p$, we write $\Sp_p$ and $\Sp_p^G$ for the categories
of $p$-complete spectra and $p$-complete genuine $G$-spectra,
respectively, equipped with the $p$-completed tensor product. A genuine spectrum is $p$-complete precisely if all the fixed point spectra are $p$-complete.
We denote inflation from the trivial group by
$\infl_e^G:\Sp\to\SpG$.

For a symmetric monoidal $\infty$-category $\mathcal V$, we write
\[
 \coAlg(\mathcal V):=\CAlg(\mathcal V^{\op})^{\op}
\]
for its category of counital $\mathbb{E}_\infty$-coalgebras. Generally all coalgebras appearing in this paper are cocommutative and counital.

\subsection*{Acknowledgments} The author would like to thank Phil Pützstück, Thorger Gei\ss{} and Maxime Ramzi for helpful discussions and comments on the proof. Phil Pützstück suggested the proof strategy that we implemented in Section \ref{sec:def_norm} using product preservation (replacing earlier Pro-isomorphism arguments that I had used). 
Special thanks to John Rognes for many discussions over the years, his amazing work on the subject and comments on a draft of this paper. I am also very grateful to Allen Yuan for discussions of coalgebras and for  sharing his unpublished work on global coalgebras and relative Segal conjectures.

The author was funded by the Deutsche Forschungsgemeinschaft (DFG, German Research
Foundation) -- Project-ID 427320536 -- SFB 1442, as well as under Germany’s Excellence Strategy
EXC2044/2–390685587, Mathematics M\"unster: Dynamics–Geometry–Structure.

\subsection*{AI Statement}I used AI to search and understand the literature, assess the existing
results, and refine the proofs and especially the exposition, in particular by
generating referee-style reports. 
Specifically, it was an AI system (Astra) which pointed out
that a key step in the argument of Bergsaker--Rognes is the
reduction to the $p$-adic vanishing of $\Psi^G(\NormG X)$
\cite[Proposition~2.2]{BR}, after which it was easy for me  to adapt the
strategy from \cite[Section~III.1]{NS} to general groups and prove
Theorem~\ref{thm:main}. For the second part, I used AI to improve
parts of the proof of Proposition~\ref{prop:sections-equivalence}. 

\part{The generalized Segal conjecture}\label{part:segal}

We now prove Theorems~\ref{Tate_diagonal} and~\ref{thm:main}.
In Section~\ref{sec:equivalent}, we reduce both statements to the
$p$-group case, where it suffices to show that $(\NormG X)^\wedge_p$
is Borel-complete. We approach this through the isotropy defect,
whose basic properties are established in
Section~\ref{sec:isotropy defect}. In Section~\ref{sec:def_norm},
we reduce its $p$-adic vanishing to the case $X=\Fp$,
established by Bergsaker--Rognes~\cite{BR}.

\section{Equivalent formulations}\label{sec:equivalent}

Recall that the Burnside ring $A(G)=\pi_0^G\Sph_G$ is the group completion of isomorphism classes of finite $G$-sets under disjoint union. It comes with the augmentation $A(G) \to \mathbb{Z}$ induced by cardinality. The augmentation ideal is the kernel of this map and is finitely generated $I(G)=(a_1,\ldots,a_r)$. For any genuine $G$-spectrum $Y$ we have
\[
 Y^\wedge_{I(G)}\simeq
 \invlim_n Y/(a_1^n,\ldots,a_r^n),
\]
where the quotient denotes iterated cofibers of multiplication by the respective elements.
This induces the $I(G)$-completion on all fixed point spectra, using the module structure over $\Sph_G^G$ by restriction. In particular $Y\to Y^\wedge_{I(G)}$ is an equivalence on underlying
spectra, since each $a_i$ restricts to zero at the trivial subgroup.
Borel completion commutes with this completion, since it preserves limits
and is exact. 
Thus every Borel-complete spectrum is $I(G)$-complete, and the Borel
comparison factors canonically through $Y^\wedge_{I(G)}$:
\[
Y \to Y^\wedge_{I(G)} \to Y^\wedge_{\Tt} \ .
\]

For bounded-below $Y$, the effect of completion on geometric fixed points is
given by \cite[Theorem 4.1 and Lemma 5.4]{Rag}:
\begin{equation}\label{eq:marks}
 \Phi^H(Y^\wedge_{I(G)})\simeq
 \begin{cases}
  \Phi^eY = Y^e &\text{if } H=e,\\
  (\Phi^HY)^\wedge_p &\text{if }H\ne e\text{ is a }p\text{-group},\\
  0&\text{if } H\text{ is not a }p\text{-group for any }p.
 \end{cases}
\end{equation}
We get similar statements for completion at the ideal $(p)$, e.g. that $p$-completion of $G$-spectra is  simply given by completing all the fixed point spectra. On geometric fixed points we find
$\Phi^H(Y^\wedge_p)\simeq(\Phi^HY)^\wedge_p$ provided that $Y$ is bounded below. 
In fact a genuine $G$-spectrum $Y$ is connective if and only if all
$\Phi^H Y$ are connective and in this case if all homotopy groups of  $\Phi^H Y$ are finite $p$-groups resp. derived $p$-complete then
the same holds for all homotopy groups of $Y^H$. Both of these facts are standard and easily follow from the usual isotropy separation sequence by induction on the group order (using that the classifying anima are of finite type). We also note that the notion of being bounded below for a $G$-spectrum is equivalent to all fixed point spectra (or all geometric fixed point spectra) being bounded below. Whether or not one asks this to be uniform through all subgroups doesn't matter since there are only finitely many subgroups.

Recall also that the proper Tate construction of a spectrum $M$ with $H$-action
is obtained by applying $\Phi^H$ to its Borel-complete genuine refinement.
It vanishes when $H$ is not a $p$-group for any prime $p$
\cite[Proposition 7.7]{Nik}. For a non-trivial $p$-group $H$, it kills
$p$-invertible spectra as well: the ambidextrous adjunction between induction and restriction gives
equivariant maps $M\to\Ind_e^H\Res_e^H M\to M$ whose composite is multiplication by
$|H|$ on underlying spectra. If $p$ acts invertibly, this composite is an
equivalence, so $M$ is a retract of an induced spectrum, which is killed by proper Tate.

We use the following identities for the HHR norm \cite[Proposition 2.1]{BR}
\begin{equation}\label{eq:norm}
 \Res_H^G\NormG X\simeq N_e^H(X^{\otimes[G:H]}),
 \qquad \Phi^H\NormG X\simeq X^{\otimes[G:H]}.
\end{equation}
In particular, $\NormG X$ is bounded below whenever $X$ is bounded below
\cite[Lemma 6.3]{BR}.

\begin{prop}\label{prop:equivalent}Theorem ~\ref{Tate_diagonal} holds for all finite $p$-groups precisely if Theorem~\ref{thm:main} holds for all finite groups\footnote{We do not claim that they are equivalent if we fix the group $G$. What we precisely show is that Theorem~\ref{thm:main} implies  Theorem~\ref{Tate_diagonal} for a fixed non-trivial finite $p$-group $G$ and is conversely implied by Theorem~\ref{Tate_diagonal} for all $p$-subgroups $H \le G$.  }.
 They are
also equivalent to the following $p$-group statement: for every prime $p$,
every finite $p$-group $G$ and every bounded-below spectrum $X$, the canonical map
\begin{equation}\label{eq:padicnorm}
 (\NormG X)^\wedge_p\longrightarrow
 \bigl((\NormG X)^\wedge_p\bigr)^\wedge_{\Tt} =  \bigl((\NormG X)^\wedge_{\Tt}\bigr)^\wedge_{p} 
\end{equation}
is an equivalence.
\end{prop}
\begin{proof}
We consider the  completion map $\NormG X\to (\NormG X)^\wedge_{\Tt}$.
For $H\le G$, restriction commutes with Borel completion, so \eqref{eq:norm}
identifies the map on geometric $H$-fixed points with
\[
 X^{\otimes[G:H]}\xrightarrow{\Delta_H}((X^{\otimes[G:H]})^{\otimes H})^{\varphi H},
\]
For $H=e$ this is the identity. If Theorem~\ref{Tate_diagonal} holds, then this is a
$p$-completion for $H$ a non-trivial $p$-group, and its target is zero
if $H$ is not a $p$-group for a prime $p$. Comparing with \eqref{eq:marks} shows
that $(\NormG X)^\wedge_{I(G)}\to (\NormG X)^\wedge_{\Tt}$ is an equivalence on all geometric
fixed points, hence an equivalence. This proves Theorem~\ref{thm:main}.

Conversely assume Theorem~\ref{thm:main} holds. Then apply $\Phi^G$ to that theorem for a non-trivial $p$-group $G$. The computation
\eqref{eq:marks} identifies the resulting map with the $p$-completion of
$X$, giving Theorem~\ref{Tate_diagonal}. This shows the equivalence between these two theorems. 

Now assume Theorem~\ref{Tate_diagonal} and let $G$ be a $p$-group.
The same calculation as above  shows that $\NormG X\to (\NormG X)^\wedge_{\Tt}$ is a mod~$p$
equivalence on all geometric fixed points, hence a genuine mod~$p$
equivalence. Completing and commuting $p$-completion with Borel completion
gives \eqref{eq:padicnorm}.

Finally, assume \eqref{eq:padicnorm} and let $G\ne e$ be a $p$-group.
The fiber of $\NormG X\to (\NormG X)^\wedge_p$ is $p$-invertible on underlying naive
$G$-spectra, so proper Tate sends this map to an equivalence. Therefore
\[
 \begin{aligned}
 (X^{\otimes G})^{\varphi G}
 \simeq\Phi^G\bigl((\NormG X)^\wedge_{\Tt}\bigr) 
 \simeq\Phi^G\Bigl(\bigl((\NormG X)^\wedge_p\bigr)^\wedge_{\Tt}\Bigr) 
 \simeq\Phi^G\bigl((\NormG X)^\wedge_p\bigr)\simeq X^\wedge_p.
 \end{aligned}
\]
The identifications respect the map from $X=\Phi^G\NormG X$, so the Tate diagonal
is the completion map. This proves Theorem~\ref{Tate_diagonal}.
\end{proof}

\section{The isotropy defect}\label{sec:isotropy defect}

In this section let $G$ be a non-trivial finite group. The (stable) isotropy separation sequence is
\[
 \Sigma^\infty_+\EP\longrightarrow\Sph_G
 \longrightarrow\Sigma^\infty\tEP
\]
where $\tEP$ is defined as the cofiber of $\EP\to *$. It is an idempotent algebra and we have
$
 \Phi^G Y\simeq(\Sigma^\infty\tEP\otimes Y)^G$ as functors $\SpG \to \Sp$. 
There is a fully faithful embedding
\[
 \Xi^G:\Sp\longrightarrow\SpG,\qquad
 \Xi^G(X)=\Sigma^\infty\tEP\otimes X,
\]
where tensoring with a spectrum $X$ uses the action of spectra on $G$-spectra (implemented through colimits and stability).  $\Phi^G$ is left adjoint to
$\Xi^G$.
The essential image of $\Xi^G$ consists precisely of
those $Y$ for which $Y^H=0$ for every proper subgroup $H \lneqslant G$. On this image, categorical
and geometric $G$-fixed points agree. 

The functor $\Xi^G$ admits a further right adjoint, described in
\cite[\S2.6]{HKS}. We recall its construction and collect the
properties needed for our proof.


\begin{defn}
We denote the right adjoint of $\Xi^G$ by $\Psi^G$ and call it the
\emph{isotropy defect}. Thus
\begin{equation*}\Psi^G(Y)=F(\tEP,Y)^G \ .
\end{equation*}
\end{defn}
\begin{remark}
$\Psi^G$  is exactly the obstruction also used by Bergsaker--Rognes \cite[Proposition~2.2 and Definition~2.3]{BR}, see  a more precise discussion below. They write $\tEP = S^{\infty\rho}$ where
$\rho=\mathbb C[G]\ominus \mathbb C$ is the reduced regular complex
representation.
\end{remark}

Since  $\Psi^G$ is a right adjoint, it is exact and preserves all limits. In particular, it commutes with
$p$-completion which will be used later.  Applying $F(-,Y)$ to isotropy separation gives a fiber sequence
\begin{equation}\label{coisotropy_sequence}
 \Xi^G\Psi^G(Y)\longrightarrow Y\longrightarrow Y^\wedge_{\Pp}
\end{equation}
which we call the coisotropy sequence.
Thus $\Psi^G(Y)=0$ if and only if $Y$ is $\Pp$-complete, hence the name `isotropy defect'. In particular,
$\Psi^G(Y^\wedge_{\Tt})=0$.

\begin{prop}\label{prop:psi}
For a proper subgroup $H\lneqslant G$, the functor $\Psi^G$
vanishes on $\Ind_H^G Z$. It preserves direct sums of genuine
spectra with a common upper bound on their homotopy groups.
\end{prop}
\begin{proof}
Induction agrees with coinduction, so
$F(\EP_+,\Ind_H^G Z)\simeq\Ind_H^G Z$, since the restriction of
$\EP$ to $H$ is contractible. The vanishing follows from
\eqref{coisotropy_sequence}.

Choose a $G$-CW model of $\tEP$ with finitely many cells in each
dimension. For $Y\in\SpG_{\le n}$ and fixed $k$, the group
$\pi_k F(\tEP,Y)^G$ is computed from a finite skeleton, whose size
depends only on $n$ and $k$. Thus $\Psi^G$ preserves direct sums
with a common upper bound.
\end{proof}

To relate the isotropy defect to the proof strategy in
\cite[Section~III.1]{NS}, consider the functor
\[
 T_G:\Sp\longrightarrow\Sp,\qquad
 T_G(X)=\Phi^G\bigl((\NormG X)^\wedge_{\Pp}\bigr).
\]
This is the `correct' generalization of the functor $T_p$ considered there. 
Applying $\Phi^G$ to \eqref{coisotropy_sequence} and using \eqref{eq:norm} we get a
natural fiber sequence
\[
 \Psi^G(\NormG X)\longrightarrow X\longrightarrow T_G(X).
\]
The proper Tate diagonal is the composite of $X\to T_G(X)$ with
the natural comparison $T_G(X)\to(X^{\otimes G})^{\varphi G}$.

The significance of these functors is the  following isotropy-separation reduction argument of
Bergsaker--Rognes~\cite[Proposition~2.2]{BR}, following
Carlsson~\cite[Theorem~A(b)]{Car} and
Caruso--May--Priddy~\cite[Lemma~1.9]{CMP}.
We recall the argument, which applies here without finite-type assumptions.
\begin{prop}\label{prop:defecteq} The $p$-group statement \eqref{eq:padicnorm} and hence Theorems~\ref{Tate_diagonal} and~\ref{thm:main} are equivalent to the assertion that  $\Psi^G( \NormG X)$ is $p$-adically trivial  for each bounded-below spectrum $X$ and each non-trivial finite $p$-group $G$. Equivalently, the map $X\to T_G(X)$ is a $p$-adic equivalence for all such $G$ and $X$.
\end{prop}

\begin{proof}
If the $p$-completion of the norm is Borel-complete, then it is also $\Pp$-complete so that by the coisotropy sequence \eqref{coisotropy_sequence} the isotropy defect vanishes $p$-adically (using that $\Psi^G$ commutes with $p$-completion).

Conversely, if $\Psi^G(\NormG X)^\wedge_p = \Psi^G((\NormG X)^\wedge_p) = 0$, then the spectrum $(\NormG X)^\wedge_p$ is $\Pp$-complete. To show that it is Borel-complete it therefore suffices to show that the map
\[
\NormG X \to (\NormG X)^\wedge_{\Tt}
\]
is a $p$-adic equivalence on all fixed point spectra for proper subgroups, i.e. after restriction to proper subgroups.
This follows by induction on $|G|$, the case $G=e$ being automatic, since
these restrictions are themselves norms by \eqref{eq:norm}.
\end{proof}

\section{The vanishing of the defect}\label{sec:def_norm}

We now prove the $p$-adic vanishing of $\Psi^G(\NormG X)$
by successive reductions, following the strategy of the proof of
\cite[Theorem~III.1.7]{NS}. Our presentation avoids the explicit use of pro-isomorphisms.
Instead, we use preservation of infinite products to pass from
$\mathbb F_p$ to arbitrary $\mathbb F_p$-vector spaces,
each of which is a retract of a product of copies of $\mathbb F_p$.
The idea of reformulating the $C_p$-proof
in this way is due to Phil Pützstück. Here we adapt it to arbitrary
finite $p$-groups. It is however not very complicated to rewrite it in a way closer to \cite{NS}  using pro-isomorphisms. We leave the details to the reader.

All arguments in this section are formal and elementary.
The (hard) result of Bergsaker--Rognes~\cite[Theorem~1.1]{BR} is invoked
only at the very last step to establish the vanishing for $\Fp$.

\begin{lemma}\label{lemma:exact}
The functor $X\mapsto\Psi^G(\NormG X)$ is exact.
\end{lemma}
\begin{proof}
For a cofiber sequence $A\to B\to C$, the tensor-product cube
gives a finite filtration of $\NormG B$, starting with $\NormG A$
and with final layer $\NormG C$. Every intermediate layer is induced
from stabilizers of nonempty proper subsets of the transitive $G$-set
$G$. These stabilizers are proper subgroups, so $\Psi^G$ kills the
intermediate layers. Thus the sequence
\[
 \Psi^G(\NormG A)\longrightarrow\Psi^G(\NormG B)
 \longrightarrow\Psi^G(\NormG C)
\]
is a cofiber sequence.
\end{proof}

\begin{lemma}\label{lem:postnikov}
For a bounded-below spectrum $X$, the natural maps
\[
 \NormG X\xrightarrow{\ \simeq\ }
 \invlim_n\NormG(\tau_{\le n}X) \qquad \text{and} \qquad   \Psi^G(\NormG X)\xrightarrow{\ \simeq\ }
 \invlim_n\Psi^G\bigl(\NormG(\tau_{\le n}X)\bigr)
\]
are equivalences. 
\end{lemma}
\begin{proof}
By \eqref{eq:norm} and the usual yoga for total fibers the connectivity of the fibers on geometric fixed points of $\NormG X \to \NormG(\tau_{\le n}X)$ tends to $\infty$, hence also on fixed points as recalled in Section~\ref{sec:equivalent}. This shows the first claim. The second follows by limit preservation of $\Psi^G$.  
\end{proof}

To prove preservation of products, we truncate the norm in genuine
$G$-spectra. Since $\Psi^G$ preserves limits, Postnikov completeness gives
\begin{equation}\label{eq:defect-postnikov}
 \Psi^G(\NormG X)\simeq
 \invlim_n\Psi^G(\tau_{\le n}\NormG X).
\end{equation}

\begin{lemma}\label{finite_dim}
Every group $\pi_k\Psi^G(\tau_{\le n}\NormG\Fp)$ is a finite $p$-group.
\end{lemma}
\begin{proof}
Finite tensor powers of $\Fp$ are connective and of finite type.
By \eqref{eq:norm} and the facts recalled in
Section~\ref{sec:equivalent}, all fixed points of $\NormG\Fp$
have finite $p$-groups as homotopy groups. The same holds for
$Y=\tau_{\le n}\NormG\Fp$.
Since $Y$ is bounded above, $\pi_k\Psi^G(Y)$ is computed from a
finite skeleton of $\tEP$, as in the proof of
Proposition~\ref{prop:psi}. Hence it is a finite $p$-group.
\end{proof}

\begin{lemma}\label{lem:prod_preservation}
The functor $V\mapsto\Psi^G(\NormG V)$ from $\Fp$-vector spaces
to spectra preserves arbitrary products.
\end{lemma}
\begin{proof}
For a family $(V_i)_{i\in I}$ of vector spaces, distribute the norm
of $\bigoplus_i V_i$ out. The constant functions $G\to I$ give
the individual norms. All other summands are induced from proper
stabilizers. Genuine Postnikov truncation commutes with direct
sums and induction, so the truncated mixed summands are still
induced from proper subgroups. They now have a common upper
bound. Proposition~\ref{prop:psi} therefore gives
\[
 \Psi^G\!\left(\tau_{\le n}\NormG\!\left(\bigoplus_i V_i\right)\right)
 \simeq\bigoplus_i\Psi^G(\tau_{\le n}\NormG V_i).
\]
Thus $V\mapsto\pi_k\Psi^G(\tau_{\le n}\NormG V)$ is infinitely
additive and hence of the form
\[
 \pi_k\Psi^G(\tau_{\le n}\NormG V)\cong
 \pi_k\Psi^G(\tau_{\le n}\NormG\Fp)\otimes_{\Fp}V.
\]
Here additivity gives the $\Fp$-vector space structure on the
coefficient, which is finite-dimensional by Lemma~\ref{finite_dim}.
Thus every functor $V\mapsto\Psi^G(\tau_{\le n}\NormG V)$ preserves
arbitrary products. The claim then follows from
\eqref{eq:defect-postnikov}, since limits commute with products.
\end{proof}

\begin{proof}[Proof of Theorem \ref{thm:main}]
By Proposition~\ref{prop:defecteq}, we need to prove that  $\Psi^G(\NormG X) / p = 0$ for each bounded-below $X$ and each non-trivial $p$-group $G$. Using the Postnikov tower and Lemmas~\ref{lemma:exact} and~\ref{lem:postnikov} reduces this to $X$ concentrated in a single degree represented by an abelian group $A$. Then we can commute the mod $p$ reduction inside and get an extension of the $p$-torsion in $A$ and the non-derived mod $p$-reduction of $A$, both of which are $\Fp$-vector spaces. 

Thus it suffices to show that $\Psi^G(\NormG V)= 0$ for an $\Fp$-vector space $V$. 
Bergsaker--Rognes' result \cite[Theorem 1.1]{BR} shows using Proposition \ref{prop:defecteq} that $\Psi^G(\NormG \Fp) = 0$.  We use Lemma \ref{lem:prod_preservation} to deduce that it vanishes on arbitrary products of $\Fp$. But then every $\Fp$-vector space $V = \bigoplus_{i \in I} \Fp$ is a retract of $\prod_{i \in I} \Fp$. Hence $\Psi^G(\NormG V)$ also vanishes. 
\end{proof}

\part{The coalgebraic Frobenius}\label{part:frobenius}

We now turn to the proof of Theorem~\ref{thm:coalgebra-frobenius-action}, i.e. the construction of the oriented circle action on the category of bounded-below, $p$-complete coalgebras.
Fix a prime $p$. Our input from Part~\ref{part:segal} is the following
consequence of the generalized Segal conjecture (Theorem~\ref{thm:main}):
the $p$-completed norm of every bounded-below spectrum is
Borel-complete for every finite $p$-group, by
Proposition~\ref{prop:equivalent} and \eqref{eq:padicnorm}.
We first use this to construct an adjunction on coalgebras.
In Section~\ref{sec:norm-systems} and Section~\ref{sec:proofnormsystem} we assemble these adjunctions into
coherent norm systems, and in Section~\ref{sec:frobenius-action}
we apply Yuan's construction to obtain the Frobenius action. 

\section{The Frobenius through genuine coalgebras}\label{sec:adjunction}

We consider cocommutative coalgebras in genuine $G$-spectra $\Sp_p^G$ and refer to them as genuine coalgebras. Such a coalgebra has genuine equivariant
comultiplication and counit maps
$A\to A\otimes A$ and $A\to\mathbbm{1}$, with the usual
$\mathbb{E}_\infty$-coherences. No `normed' or `ultracommutative' structure
is included in this definition. In particular, maps into indexed norms
are not part of the data. 

In this and the following sections, tensor products, inflation
$\infl_e^G$, norms $N_H^G$, and geometric
fixed points are taken in the $p$-complete sense.
Thus $N_H^G Y$ means  $(N_H^G Y)^\wedge_p$ and $\Phi^HY$ means $(\Phi^HY)^\wedge_p$.
We will also  consider genuine residual actions on geometric fixed points. For example  $\Phi^H$ refines to a functor from genuine $(G\times H)$-spectra to
genuine $G$-spectra.

Norms, inflation, and geometric fixed points preserve mod~$p$
equivalences. Hence the usual norm and base-change identities
\cite[\S9.2]{BH} descend to the $p$-complete setting. 
These  functors are also strong symmetric monoidal and hence induce
functors on coalgebras. The tensor product of coalgebras is their
categorical product, and the tensor unit is terminal, see e.g.
\cite[\S2]{BB}. In particular, the induced functors preserve
finite products.

Let $G$ be a finite $p$-group, and write $U$ for the forgetful functor from
genuine spectra to spectra with $G$-action, i.e. functors from $BG$ to $\Sp_p$. 
\begin{lemma}\label{lem:borel-comparison}
Let $B\in\coAlg(\Sp_p^G)$ be such
that the underlying spectrum of $B^{\otimes n}$ is Borel-complete
for every $n\ge0$\footnote{For $n=0$ this is the $p$-complete genuine sphere which is Borel-complete under our assumptions and conventions.}.
Then, for every $A\in\coAlg(\Sp_p^G)$,
restriction induces a natural equivalence
\begin{equation*}
 \Map_{\coAlg(\Sp_p^G)}(A,B)
 \simeq
 \Map_{\Fun(BG,\coAlg_p)}(UA,UB).
\end{equation*}
In particular, restriction is fully faithful on the full subcategory
of coalgebras satisfying this condition.
\end{lemma}
\begin{proof}
Using the operadic description, we can describe the mapping anima between $\mathbb{E}_\infty$-coalgebras as limits of mapping anima between tensor powers.
All target tensor powers are Borel-complete, so restriction
identifies every term with its counterpart in $\Fun(BG,\Sp_p)$.
Taking the operadic limits gives the asserted equivalence.
The case $n=0$ includes the compatibility with the counit.
\end{proof}

We continue to let $G$ be a finite $p$-group.

\begin{lemma}\label{lem:lift} 
For any $A\in\coAlg(\Sp_p^G)$ and any bounded-below $p$-complete coalgebra
$B$, there is a natural equivalence
\[
 \Map_{\coAlg(\Sp_p^G)}(A,N_e^GB)
 \simeq
 \Map_{\coAlg_p}(\Res_e^G A,B).
\]
In particular, restriction is left adjoint to the norm on
bounded-below coalgebras. 
%
\end{lemma}
\begin{proof}
Every tensor power $B^{\otimes n}$ is bounded below.
Thus $(N_e^GB)^{\otimes n} = N_e^G(B^{\otimes n})$
is Borel-complete by \eqref{eq:padicnorm}. Lemma~\ref{lem:borel-comparison} therefore identifies mapping
anima from any $A$ into $N_e^GB$ with the counterpart after applying $U$.

Using that the tensor product in counital cocommutative coalgebras is the  categorical
product we see that the underlying naive coalgebra $B^{\otimes G}$, with
its permutation action, is coinduced. This identifies the remaining mapping anima
with $\Map_{\coAlg_p}(\Res_e^G A,B)$.
This proves the stated equivalence of mapping anima, naturally in both variables.
Since restriction, inflation and the norm preserve bounded-below objects, this also gives the adjunction.
%
\end{proof}

If $\eta$ denotes the unit of the previous adjunction, we get for every bounded-below, $p$-complete coalgebra $C$ a canonical map 
\[
 \bar{\Delta}_G:=\eta_{\infl_e^G C}: \quad 
 \infl_e^G C\longrightarrow
 N_e^G\Res_e^G\infl_e^G C
= N_e^G C \ .
\]
Concretely under the adjunction between inflation and fixed points this corresponds to a map of spectra $C \to  (N_e^GC)^G$, which by completeness of the target is simply the map induced by comultiplication $\Delta^G: C \to (C^{\otimes G})^{hG}$. 

\begin{defn}\label{def:vector-frobenius}
For a finite-dimensional $\Fp$-vector space $V$ viewed as an
elementary abelian $p$-group we
define
\[
\varphi_V:=\Phi^V  \bar{\Delta}_V:\quad \id_{\coAlg_p^{\mathrm{bb}}}\to\id_{\coAlg_p^{\mathrm{bb}}}
\]
using
$\Phi^V\infl_e^V= \id$ and $\Phi^VN_e^V = \id$. We also define $\varphi_0:=\id$.
\end{defn}
For $V\ne0$, comparison with Borel completion identifies its underlying map as
\begin{equation}\label{eq:frob}
 C\xrightarrow{\Delta^V} (C^{\otimes V})^{hV}
 \xrightarrow{\mathrm{can}}(C^{\otimes V})^{\varphi V}
 \xrightarrow{\Delta_V^{-1}}C.
\end{equation}
Here $\varphi V$ denotes proper Tate fixed points and $\Delta_V$ is
the proper Tate diagonal, which is an equivalence here by
Theorem~\ref{Tate_diagonal}. For $V=\Fp$, the proper Tate construction agrees with the ordinary Tate construction, so
$\varphi_{\Fp}=\varphi_p$ is precisely the coalgebraic Frobenius of
\eqref{eq:coalgebra-frobenius}. 

\begin{example}
For $V=\Fp\oplus\Fp$, the endomorphism $\varphi_V$ is the composite
of the coalgebraic Frobenius with itself:
\[
 \varphi_{\Fp\oplus\Fp}\simeq\varphi_p\circ\varphi_p.
\]
This follows since the $V$-fold comultiplication is obtained by applying the
$p$-fold comultiplication twice. Using Lemma~\ref{lem:lift} gives the corresponding factorization of its
genuine lift. Taking geometric fixed points successively in the two
$\Fp$-directions thus shows the assertion, see also 
\cite[Examples 2.11 and 2.12]{Yuan}.
\end{example}

The idea of our proof of 
Theorem~\ref{thm:coalgebra-frobenius-action} (inspired by Yuan's work) 
is to first make such
identifications 
$ \varphi_{\Fp\oplus\Fp}\simeq\varphi_p\circ\varphi_p$
coherently compatible with maps and direct
sums of vector spaces, thus isolating the ways in which the coalgebraic Frobenius ``commutes with itself''.
This yields a monoidal action of a
category built from finite-dimensional $\Fp$-vector spaces
on $\coAlg_p^{\mathrm{bb}}$, precisely the monoidal category $Q\Vect$ discussed at the beginning of Section \ref{sec:frobenius-action}.
Then, following Yuan, we use his partial $K$-theory analogue
of Quillen's computation for $\Fp$
\cite[Theorem 5.2]{Yuan}, together with $p$-completeness,
to descend this action to the oriented circle, i.e. to obtain
an action of $B\mathbb N$.

\section{Coalgebraic norm systems}\label{sec:norm-systems}

We now recall the correspondence categories from Yuan's paper~\cite{Yuan},
with some modifications that we need for the present paper. Generally we review what we need from the work of Yuan, but assume familiarity with the ideas and methods.

We work with finite groupoids with elementary abelian $p$-isotropy,
i.e. finite disjoint unions (possibly empty) of $BV$ for
finite-dimensional $\Fp$-vector spaces $V$, including $V=0$.
Here $BV$ denotes the classifying anima of the additive group of $V$.
We consider those as a $(2,1)$-category, defined as a full subcategory of the category of anima. 
A \emph{finite covering}
is a map between such groupoids whose fibers are finite discrete anima. For example,
$*\to BV$ is a finite covering, but $BV\to *$ is not (unless $V=0$).

\begin{defn}\label{def:corr}
The category $\Corr$ has these groupoids as objects and spans
\[
 a=(X\xleftarrow{f}M\xrightarrow{g}Y)
\]
with $g$ a finite covering
as morphisms. Composition is by  pullback. Formally this is a  subcategory of the span category of anima (see e.g.\ \cite[\S2]{CLR}). 
\end{defn}

\begin{defn}\label{def:corr-plus}
The category $\Corr^{+}$ has the same objects, but its morphisms are
spans for which the joint map $(f,g):M\to X\times Y$ is a finite
covering. To compose $X\leftarrow M\to Y$ and $Y\leftarrow N\to Z$,
factor the natural map as
\[
 M\times_Y N\longrightarrow P\longrightarrow X\times Z,
\]
where the first map has connected fibers and the second
is a finite covering. The composite is $X\leftarrow P\to Z$.
Thus $P$ is the relative $0$-truncation of $M\times_Y N$ over
$X\times Z$. This composition restricts to the usual composition
on $\Corr$, giving a wide subcategory $\Corr\subseteq\Corr^{+}$.
A precise construction of this $(2,1)$-category is given in
\cite[Definition~2.13 and Convention~2.17]{Yuan}, from which we obtain
$\Corr^{+}$ by restricting to the groupoids considered here.
\end{defn}

\begin{defn}\label{def:corr-pointed}
The category $\Corr^{*}$ has pointed groupoids $(X,x)$ as objects
and spans in $\Corr$ with a chosen point in $M$ and point-preserving
maps to $X$ and $Y$ as morphisms. Composition is by
pullback, i.e. this is a subcategory of the span category of pointed anima.
\end{defn}

These are $(2,1)$-categories but we regard them as $\infty$-categories here.
We call a morphism in $\Corr$ or $\Corr^{*}$ a backward morphism if  $g=\id$ and a
forward morphism if $f=\id$.
We have functors
\[
 \Corr^{*}\longrightarrow\Corr\subseteq\Corr^{+}.
\]
Cartesian products of groupoids (which are not categorical products in the span categories) give symmetric monoidal structures
on these categories so that these functors are symmetric monoidal.

\begin{remark} Our $\Corr$ and $\Corr^{+}$ are the full
subcategories of Yuan's $\mathrm{Glo}$ and $\mathrm{Glo}^{+}$,
respectively, on groupoids with elementary abelian $p$-isotropy
\cite[Variant 2.5 and Definition 2.20]{Yuan}. $\Corr^*$ has not explicitly appeared in Yuan's work. We warn the reader that $\Corr^*$ is not the category of pointed objects in $\Corr$ as the latter is already pointed. 
\end{remark}

\begin{construction}\label{cons:qplus}
We now assemble genuine coalgebras into a category over $\Corr^{*}$.
For $X\simeq\coprod_i BG_i$, write
$\Sp_p^X:=\prod_i\Sp_p^{G_i}$ for the category of $p$-complete
genuine spectra over $X$. The proof of \cite[Theorem 2.18]{Yuan}
assembles genuine spectra into a coCartesian fibration over $\Corr^{+}$ with symmetric monoidal fibers and strong symmetric monoidal transition
functors. Localizing this system at the mod $p$-equivalences as
in Section~\ref{sec:adjunction} and applying $\coAlg$ fiberwise gives
a coCartesian fibration
\[
 q^{+}:\mathcal E^{+}\longrightarrow\Corr^{+},
 \qquad \mathcal E^{+}_X = \coAlg(\Sp_p^X).
\]
The fiber transport along a span
$(X\xleftarrow{f}M\xrightarrow{g}Y)$ is
\[
 g_\otimes f^*:\coAlg(\Sp_p^X)\longrightarrow\coAlg(\Sp_p^Y),
\]
where $f^*$ is restriction and $g_\otimes$ is the multiplicative pushforward along $g$ (which specializes to norms for finite coverings, in particular for subgroup inclusions,
and to geometric fixed points for quotient homomorphisms, see \cite[\S2.1]{Yuan}). Pullback along the functor $\Corr^{*}\to\Corr^{+}$ yields a coCartesian fibration
\begin{equation}\label{eq:coalgebra-fibration}
 q:  \mathcal E\longrightarrow\Corr^{*},
 \qquad \mathcal E_{(X,x)} = \coAlg(\Sp_p^X).
\end{equation}
Thus an object of $\mathcal E$ is a pointed groupoid $(X,x)$ together
with a genuine $X$-coalgebra and a morphism is a span of pointed groupoids $(f,g)$ together with 
a coalgebra map $g_\otimes f^*A\to B$. Here $g$ is a finite covering, so that $g_\otimes$ is indeed the norm for a covering. 
\end{construction}

\begin{defn}\label{def:coalgebraic-norm-systems}
A \emph{coalgebraic norm system}, or simply a \emph{norm system},
is a section of $q$ (as in 
\eqref{eq:coalgebra-fibration}) which is Cartesian
over backward morphisms. We write
$\NormSys_p$ for the full $\infty$-category of
these sections. 
\end{defn}

\begin{remark}
\begin{itemize}
\item
Informally, a norm system packages coalgebras for all
elementary abelian $p$-groups  into one object: it assigns a genuine $V$-coalgebra $s(BV)$ to each such
group $V$, together with maps relating these coalgebras under
restriction and norms. 
Proposition~\ref{prop:sections-equivalence} below shows that, when
$C=s(*)$ is bounded below, this entire coherent system is determined
by $C$: it is a coherent system of norm powers, with
$s(BV)\simeq N_e^VC$. Therefore we think of these systems as controlling the coherence and functoriality of norms of elementary abelian $p$-groups. 
\item

Unlike Yuan's global algebras \cite[Definition~3.19]{Yuan},
which inspired our definition, our norm systems are Cartesian
rather than coCartesian over backward morphisms.
Pointings are essential here: over $\Corr$,
Cartesianness over the backward spans
$X\leftarrow\varnothing\to\varnothing$ would force every $s(X)$
to be terminal, since the fiber over $\varnothing$ is terminal.
The chosen points exclude these spans and provide the finite
coverings $x:*\to X$ used to form the norm powers below.
\end{itemize}
\end{remark}

Explicitly, for a pointed map $f:X\to Y$, Cartesianness over the backward morphism means that
the section map $f^*s(Y)\to s(X)$ induces an equivalence
\[
 \Map_{\coAlg(\Sp_p^Y)}(A,s(Y))\simeq\Map_{\coAlg(\Sp_p^X)}(f^*A,s(X))
\]
for every $A$. Note that  $q$ is not a Cartesian fibration. 
In the following we also use forward coCartesianness over forward morphisms associated with a finite pointed covering $g:X\to Y$. This then concretely means that
$g_\otimes s(X)\to s(Y)$ is an equivalence. We will see below that this is automatic for a norm system, provided that $s(*)$ is bounded below.

For a pointed groupoid $(X,x)$, the point $x:*\to X$ is a finite
covering. Write $N_x:=x_\otimes$ for the norm along $x$.
On the marked component $BG$, $N_xC$ is $N_e^GC$,
and on every other component it is the terminal
coalgebra $\mathbbm{1}$. Lemma~\ref{lem:lift} gives the equivalence
\begin{equation}\label{eq:point-adjunction}
 \Map_{\coAlg(\Sp_p^X)}(A,N_xC)
 \simeq \Map_{\coAlg_p}(x^*A,C)
\end{equation}
for every genuine $X$-coalgebra $A$ and bounded-below $C$.
Its counit is projection to the distinguished factor of
$x^*N_xC\simeq C^{\otimes\Omega X}$. 

The main result about norm systems is the following:

\begin{prop}\label{prop:sections-equivalence}
Let $\NormSys_p^{\mathrm{bb}}$ denote the full
subcategory of norm systems $s$ for which $s(*)$ is bounded
below. Evaluation at the point induces an equivalence
\[
 \mathrm{ev}_*:\NormSys_p^{\mathrm{bb}}
 \xrightarrow{\simeq}\coAlg_p^{\mathrm{bb}}.
\]
Its inverse sends $C$ to a section with value $N_xC$ at $(X,x)$.
Every such section is also coCartesian over forward morphisms.
\end{prop}

\section{Proof of Proposition \ref{prop:sections-equivalence}}\label{sec:proofnormsystem}

The proof has three steps. First, we prove a general extension statement
for sections over span categories (Lemmas \ref{lem:sections-yoneda} and \ref{lem:sections-span-extension}), then apply it to coalgebras in spectra
with group action, i.e. `naive' coalgebras (Lemma \ref{lem:naive-sections-equivalence}). For bounded-below coalgebras, we then use the
generalized Segal conjecture in the form \eqref{eq:padicnorm} and
Lemma~\ref{lem:borel-comparison} to lift this system uniquely to
genuine coalgebras finishing the proof of Proposition \ref{prop:sections-equivalence}.

We first need a general description of Cartesian sections in a coCartesian fibration. 

\begin{lemma}\label{lem:sections-yoneda}
Let $q:\mathcal E\to B$ be a coCartesian fibration classified by the functor 
$F:B\to\operatorname{Cat}_\infty$. Then the category of sections of $q$ that  are $q$-Cartesian over a chosen
class of morphisms in $B$ is equivalent to  the full subcategory of
$\Fun(B,\operatorname{Cat}_\infty)_{/F}$ consisting of those natural transformations  whose components  are representable right fibrations
and whose naturality squares over the chosen morphisms in $B$ are pullbacks.
The equivalence is implemented by the assignment 
\[
 s\longmapsto(F_s\to F),\qquad F_s(b)=F(b)_{/s(b)}.
\]
\end{lemma}

\begin{proof}
Let $\mathrm{RSl}\subseteq\Fun(\Delta^1,\operatorname{Cat}_\infty)$
be the full subcategory of representable right fibrations.
Then it is shown in \cite[Corollaries~6.4--6.5]{CH} that evaluation at the target
$\mathrm{RSl}\to\operatorname{Cat}_\infty$ is the universal
coCartesian fibration. Thus we have 
$\mathcal E\simeq B\times_{\operatorname{Cat}_\infty}\mathrm{RSl}$.
We apply $\Fun(B,-)$ and take the
fiber over $\id_B$ to get
\[
 \Fun_{/B}(B,\mathcal E)\simeq
 \Fun(B,\mathrm{RSl})
 \times_{\Fun(B,\operatorname{Cat}_\infty)}\{F\}.
\]
The left side is the category of sections. The right 
side is the full subcategory of
$\Fun(B,\operatorname{Cat}_\infty)_{/F}$ whose components are
representable right fibrations. This shows the stated equivalence without an additional class in the base.

For a given class in the base, the naturality square of $F_s\to F$ at a morphism $f: b \to c$ in this class is a pullback if and only if the induced map
\[
 F(b)_{/s(b)}\longrightarrow
 F(b)\times_{F(c)}F(c)_{/s(c)}
\]
is an equivalence. These are total spaces of right fibrations over $F(b)$, so being an equivalence
can be tested  fiberwise over $F(b)$. This boils down to checking whether the map
\[
 \Map_{F(b)}(z,s(b))\longrightarrow
 \Map_{F(c)}(F(f)z,s(c))
\]
is an equivalence for every $z\in F(b)$. But this is clearly equivalent to $q$-Cartesianness of $s(f)$.\end{proof}

We now use methods for extending functors to span categories
developed in the theory of six-functor formalisms, specifically
\cite[Theorem~3.3]{DK}.
To this end let $D$ be a category with pullbacks and $P$ a class of morphisms
containing all equivalences and closed under pullbacks and composition.
Suppose that $g,gf\in P$ implies $f\in P$, and that each morphism
in $P$ is $n$-truncated for some $n$ that is allowed to depend on the morphism.
Let $F:D^{\op}\to\operatorname{Cat}_\infty$ be a functor. We write
$f^*=F(f)$ and assume that $f^*$ has a right adjoint $f_*$ for $f\in P$
and that these adjoints satisfy Beck--Chevalley for pullback squares.
Let
\[
 q:\mathcal E\longrightarrow\Span_{\mathrm{all},P}(D)
\]
be the coCartesian fibration associated to the span extension $\overline{F}$ of $F$,
with transport $g_*f^*$ along the span $X\xleftarrow{f}M\xrightarrow{g}Y$ (using \cite[Section~11]{Barwick} or \cite[Corollary~2.6(3)]{CLR}).

\begin{lemma}\label{lem:sections-span-extension}
Restriction to $D^{\op}$ induces an equivalence between the full
subcategories of sections of $q$ and $q|_{D^{\op}}$ that are,
in both cases, Cartesian over backward morphisms associated
to maps in $P$.
Every such section of $q$ is also coCartesian over all forward
morphisms.

This equivalence restricts to an equivalence between the full
subcategories of sections Cartesian over all backward morphisms.
\end{lemma}

\begin{proof} We want to apply the first equivalence of \cite[Theorem~3.3]{DK} with
$E=P$ and $I$ the equivalences in $D$. The assumptions on $P$ make
this a Nagata setup in the sense of \cite[Definition~2.1]{DK}. Thus the cited theorem gives an equivalence between   $\Fun'(D^{\op},\operatorname{Cat}_\infty)$ and $\Fun''(\Span_{\mathrm{all},P}(D),\operatorname{Cat}_\infty)$.
These are certain subcategories of the ambient functor categories that we will describe now:
\begin{itemize}
\item
For $\Fun'$ the objects are functors for which
$f^*$ has a right adjoint $f_*$ for every $f\in P$, satisfying
Beck--Chevalley. The morphisms are natural transformations
whose naturality squares over maps in $P$ are right adjointable
\cite[Definition~3.1]{DK}.
\item
For $\Fun''$ write $f_!$ for forward transport.
For $f:X\to Y$ in $P$, backward transport $f^*$ must admit a
right adjoint $f_*$. Let $p_1,p_2:X\times_Y X\to X$ be the projections.
The canonical mate
\[
 f_!p_{2,*}\longrightarrow f_*p_{1,!}
\]
of the base-change equivalence
$f^*f_!\simeq p_{1,!}p_2^*$ must be an equivalence.
We require this also after arbitrary base change and recursively
for diagonals. This is stable (Pr-II) in
\cite[Definition~2.6]{DK}.
These functors form a full subcategory, so all natural transformations
between them are allowed \cite[Definition~2.15]{DK}.
\end{itemize}
Note that $F\in\Fun'$ by assumption and $\overline F\in\Fun''$ corresponds to it 
by the equivalence and so we get an equivalence on slices
\begin{equation}\label{eq:slices}
 \Fun'(D^{\op},\operatorname{Cat}_\infty)_{/F}
 \simeq 
 \Fun''(\Span_{\mathrm{all},P}(D),\operatorname{Cat}_\infty)_{/\overline F}.
\end{equation}
By Lemma~\ref{lem:sections-yoneda}, the categories of sections of $q|_{D^{\op}}$ and $q$
are equivalent to full subcategories of
$\Fun(D^{\op},\operatorname{Cat}_\infty)_{/F}$
and $\Fun(\Span_{\mathrm{all},P}(D),\operatorname{Cat}_\infty)_{/\overline F}$, respectively.
We claim that the full subcategories encoding sections Cartesian
over backward morphisms associated to maps in $P$ lie in $
 \Fun'(D^{\op},\operatorname{Cat}_\infty)_{/F}$ and $\Fun''(\Span_{\mathrm{all},P}(D),\operatorname{Cat}_\infty)_{/\overline F}$
respectively, and correspond under the equivalence \eqref{eq:slices}. This then proves the first claim.

\begin{itemize}
\item
We first consider the $\Fun'(D^{\op},\operatorname{Cat}_\infty)$-side and show that backward
Cartesianness of $s$ is equivalent to right adjointability of
the naturality squares of $F_s\to F$ over maps in $P$.
For $f:X\to Y$ in $P$,
a section $s$ is Cartesian over the associated backward morphism
if and only if the map $s(Y)\to f_*s(X)$, adjoint to the structure
map $f^*s(Y)\to s(X)$, is an equivalence. In this case, the induced functor
$F(Y)_{/s(Y)}\to F(X)_{/s(X)}$
has right adjoint
\[
 (B\to s(X))\longmapsto
 \bigl(f_*B\to f_*s(X)\simeq s(Y)\bigr).
\]
After forgetting the maps to $s(X)$ and $s(Y)$, this adjunction
has the same unit and counit as $f^*\dashv f_*$.
Hence the naturality square of $F_s\to F$ is right adjointable,
as required for a morphism in $\Fun'$.
Conversely, if the forgetful square is
right adjointable, evaluating the relevant mate at the terminal
object $\id_{s(X)}$ gives an equivalence $s(Y)\to f_*s(X)$.

The forgetful functors $F_s(X)=F(X)_{/s(X)}\to F(X)$ send the
Beck--Chevalley transformations for $F_s$ to those of $F$.
They are conservative, so $F_s$ lies in $\Fun'$.
A map between sections Cartesian over backward morphisms from $P$
induces postcomposition functors between the slices.
The corresponding mates become identities after
forgetting, so every such map satisfies the morphism condition
in $\Fun'$.
Together, this proves the claimed inclusion of the category
encoding sections Cartesian over backward morphisms from $P$
into $\Fun'(D^{\op},\operatorname{Cat}_\infty)_{/F}$.

\item
On the $\Fun(\Span_{\mathrm{all},P}(D),\operatorname{Cat}_\infty)_{/\overline F}$
side, let $s$ be a section which is Cartesian over the backward
morphisms from $P$.
The same adjunctions and forgetful functors as above send the relevant
(Pr-II) mate to the canonical comparison isomorphism
$f_*p_{2,*}\to f_*p_{1,*}$ (induced by 
$fp_1=fp_2$).
By conservativity of the forgetful functors, the (Pr-II)
mate for $F_s$ is therefore also an equivalence.
The same argument applies after base change and to iterated
diagonals, which are also in $P$ by pullback closure and the 
cancellation property. This recursively proves stable (Pr-II) and therefore that $F_s$ belongs to $\Fun''$.
Since $\Fun''$ is full we do not need to check anything on morphisms. 
\end{itemize}

We have identified the relevant subcategories as full
subcategories of the two sides of \eqref{eq:slices}.
The inclusion $D^{\op}\to\Span_{\mathrm{all},P}(D)$ is the identity
on objects, so restriction preserves and detects objectwise
representability. The characterization of backward
Cartesianness above then shows that \eqref{eq:slices} induces
the equivalence of the lemma. Cartesianness over all backward morphisms gives a further
restriction of this equivalence, since those morphisms are already
in $D^{\op}$.

Finally, the span extension identifies forward transport on the
slices with the right adjoint of backward transport.
It therefore preserves terminal objects.
Since it sends $\id_{s(X)}$ to the structure map
$f_*s(X)\to s(Y)$, this map is an equivalence, proving forward
coCartesianness.
\end{proof}

We now apply the previous lemma to `naive' coalgebras. To this end let $D$ be the category of
pointed finite groupoids with elementary abelian $p$-isotropy and
$P$ the class of finite coverings, so that we have $\Corr^*=\Span_{\mathrm{all},P}(D)$.
We consider the contravariant functor 
\[
 X \mapsto F(X)=\coAlg(\Fun(X,\Sp_p)) \simeq \Fun(X,\coAlg_p),
\]
where $\Fun(X,\Sp_p)$ carries the pointwise $p$-complete tensor product.
The functoriality is given by pullback: for $f:X\to Y$ it is $f^*:F(Y)\to F(X)$.
For $f\in P$, i.e. a finite covering, the right adjoint $f_*$ on coalgebras is given
by taking products in $\coAlg_p$, hence by tensor
products. It satisfies Beck--Chevalley.
In the associated coCartesian fibration
$q_{\mathrm{naive}}:\mathcal E_{\mathrm{naive}}\to\Corr^*$ the transport along a span
 $X\xleftarrow{f}M\xrightarrow{g}Y$ with $g \in P$ is given by $g_*f^*$.

\begin{lemma}\label{lem:naive-sections-equivalence}
Evaluation at $*$ gives an equivalence from the full category
of backward Cartesian sections of $q_{\mathrm{naive}}$ to $\coAlg_p$.
The inverse sends a coalgebra $C$ to a section $\sigma_C$ with
$\sigma_C(X,x)=x_*C$. The section $\sigma_C$ is automatically
coCartesian over forward morphisms.
\end{lemma}

\begin{proof}
Finite coverings satisfy the hypotheses of
Lemma~\ref{lem:sections-span-extension}, which reduces the claim
to Cartesian sections of $q_{\mathrm{naive}}|_{D^{\op}}$
and gives forward coCartesianness.
The point $*$ is terminal in $D^{\op}$ and every map
$x:*\to X$ is a finite covering. Therefore every $C\in\coAlg_p$
admits a relative right Kan extension from $*$ with values $x_*C$.
Evaluation at $*$ identifies the category of these extensions
with $\coAlg_p$
\cite[Proposition~4.3.2.15 and Remark~4.3.2.18]{HTT}.
By cancellation of Cartesian morphisms, these extensions are
precisely the Cartesian sections, proving the claim.
\end{proof}

\begin{proof}[Proof of Proposition~\ref{prop:sections-equivalence}]
We continue to let $\mathcal E$ be the total space of the genuine coalgebra fibration over $\Corr^{*}$, see \eqref{eq:coalgebra-fibration}.
Let
$\mathcal E_{\mathrm B}\subseteq\mathcal E$ be the full
subcategory on objects whose underlying genuine spectra and all their
finite tensor powers, including the zeroth power, are
Borel-complete. Restriction to spectra with action induces a functor
\[
 r:\mathcal E_{\mathrm B}\longrightarrow\mathcal E_{\mathrm{naive}}
\]
 over $\Corr^*$. 
This functor is fully faithful on total categories
since over a span $(X\xleftarrow{f}M\xrightarrow{g}Y)$,
Lemma~\ref{lem:borel-comparison} gives
\[
 \Map_{\coAlg(\Sp_p^Y)}(g_\otimes f^*A,B)
 \simeq
 \Map_{\coAlg(\Fun(Y,\Sp_p))}(g_*f^*rA,rB).
\]
Here we have used  that Lemma~\ref{lem:borel-comparison} 
only requires
Borel-completeness of the target tensor powers (as $g_\otimes f^*A$ is generally not Borel complete) and that the underlying spectrum of the genuine norm $g_\otimes$ is the tensor product, which is the product in $\coAlg_p$ and thus $g_*$. 

For bounded-below $C$, the generalized Segal conjecture \eqref{eq:padicnorm} applied on each
component of $X$ implies that
$
 (N_xC)^{\otimes n}\simeq N_x(C^{\otimes n})
$
is Borel-complete for each $n\ge0$. Moreover we have $r(N_xC)\simeq x_*C$.
Thus  the naive section $\sigma_C$ in fact induces a genuine section.  \\

For a pointed map $f:(X,x)\to(Y,y)$, the lifted backward map
$f^*N_yC\to N_xC$ is adjoint to the counit
$x^*f^*N_yC\simeq y^*N_yC\to C$, as follows by restriction
and full faithfulness of $r$.
By \eqref{eq:point-adjunction}, it therefore induces
\[
 \Map_{\coAlg(\Sp_p^Y)}(A,N_yC)
 \simeq\Map_{\coAlg_p}(y^*A,C)
 \simeq\Map_{\coAlg(\Sp_p^X)}(f^*A,N_xC)
\]
for every genuine $Y$-coalgebra $A$. This proves backward
Cartesianness.

Conversely, if $s$ is a backward Cartesian genuine section with
$C=s(*)$ bounded below, its backward point inclusion maps identify
$s(X)\simeq N_xC$. These identifications carry the backward
point inclusion maps to the point-inclusion adjunction counits, whose restrictions
are the usual coinduction counits. Hence $s$
takes values in $\mathcal E_{\mathrm B}$, and $r(s)$ is
Cartesian over backward point maps. Every backward map followed
by the backward point map is again a backward point map, so
cancellation shows that $r(s)$ is backward Cartesian everywhere.

Since a fully faithful functor induces fully faithful functors
on diagram categories, this identifies the genuine sections in
question with the naive backward Cartesian sections whose value
at $*$ is bounded below. Lemma~\ref{lem:naive-sections-equivalence}
therefore proves the claimed evaluation equivalence.

Finally, for a pointed finite covering $g:(X,x)\to(Y,y)$,
transitivity of norms gives $g_\otimes N_xC\simeq N_yC$.
Thus the source and target of the section's structure map belong to
$\mathcal E_{\mathrm B}$. This map is an equivalence naively by
forward coCartesianness, hence genuinely by full faithfulness.
This proves  forward coCartesianness as well.
\end{proof}

\section{The oriented circle action}\label{sec:frobenius-action}

We now finally want to construct the $B\mathbb{N}$ action on coalgebras. We first construct an action of another category, namely $Q\Vect$. This is  Quillen's $Q$-construction applied to the exact category $\Vect$ of finite-dimensional $\Fp$-vector spaces with injections and surjections. Concretely the objects of $Q\Vect$ are finite-dimensional vector spaces and  
the morphisms  are spans
$U\twoheadleftarrow W\hookrightarrow V$ and the $2$-morphisms are given by equivalences of spans. This is a priori a $(2,1)$-category, but it is equivalent to a 1-category. It is (symmetric) monoidal with respect to direct sum. 

We have the following analogue of \cite[Theorem~3.23]{Yuan}, with the coalgebra
fibration of Construction~\ref{cons:qplus} in place of the fibration
of genuine spectra.

\begin{prop}\label{prop:action}
There is a monoidal action of $Q\Vect$ on the category of sections of $\mathcal E \to \Corr^*$
given pointwise on objects by
\begin{equation*}
 (V \star s)(X)=\Phi^V s(X\times BV).
\end{equation*}
Here $\Phi^V$ denotes geometric fixed points relative to
$X\times BV\to X$, retaining the residual action on each component.
\end{prop}
\begin{proof}
The pointed groupoid construction at the end of the proof of
\cite[Proposition~3.24]{Yuan} restricts to a monoidal functor
$Q\Vect\to\Corr^*$. Taking products gives endofunctors
$(X,x)\mapsto(X\times BV,(x,*))$ of $\Corr^*$, with a morphism
$U\twoheadleftarrow W\hookrightarrow V$ acting by the pointed span
\[
 X\times BU\longleftarrow X\times BW\longrightarrow X\times BV.
\]
After forgetting basepoints this is Yuan's construction.
The augmentations give natural transformations
from these endofunctors composed with $i:\Corr^*\to\Corr^+$ to $i$,
with components
\[
 X\times BV\xleftarrow{\id}X\times BV\xrightarrow{\mathrm{pr}_X}X.
\]
We now apply \cite[Proposition~A.10]{Yuan}, which works for any
coCartesian fibration, although it is only stated for a subcategory inclusion. The proof in particular applies to our functor $i$: it uses only precomposition
and pointwise coCartesian transport along the augmentations.
Applied to $q^{+}$, this gives a right action on all sections
of its pullback $q$. Using the symmetry of $Q\Vect$, we regard
it as a left action. Transport along the augmentation is
given by geometric fixed points on coalgebras, which gives
the stated formula.
\end{proof}

\begin{remark}
This construction is where we use the genuine equivariant nature of the norm systems and why we cannot work  in the Borel world.
Key for the action are geometric fixed points, which
are strong symmetric monoidal.
The corresponding proper Tate construction on Borel spectra is only
lax symmetric monoidal and therefore does not by itself induce
the required functors on coalgebras. Thus the naive fibration used in
Lemma~\ref{lem:naive-sections-equivalence} does not directly give this action.
Compare \cite[Remark 3.36]{Yuan} for a similar phenomenon in the world of algebras used in Yuan's work. 
\end{remark}

\begin{prop}\label{prop:Q}
The action of Proposition \ref{prop:action} preserves the full subcategory
$\NormSys_p^{\mathrm{bb}}$ of bounded-below norm systems.
Via Proposition~\ref{prop:sections-equivalence}, it induces an action of $Q\Vect$  on $\coAlg_p^{\mathrm{bb}}$. On objects we have 
\[
V \star C=\Phi^V N_e^V C \simeq C
\]
Under these canonical identifications, the map induced by the span 
$
U\twoheadleftarrow W\hookrightarrow V$ is 
$\varphi_{\ker(W\to U)},
$
where $\varphi_K$ is the endomorphism of Definition \ref{def:vector-frobenius}.
In particular, the action inverts the forward morphisms.
\end{prop}
\begin{proof}
Let $s\in\NormSys_p^{\mathrm{bb}}$ and put $C=s(*)$.
Transitivity of norms and compatibility of geometric fixed points with norms give
\[
 (V\star s)(X) \simeq N_x\bigl(\Phi^V N_e^VC\bigr)
\]
where  the components away from the marked component are terminal.
In particular, its value at $*$ is $\Phi^V N_e^VC\simeq C$,
so it is again bounded below.

For backward Cartesianness it suffices to check the point maps
$x:*\to X$ since  a backward map followed by the backward pointing
map is the backward pointing map for its source, and Cartesian edges
have the relevant cancellation property. On the marked
component $BG$ the relevant map before applying $\Phi^V$ is
\[
 \Res_V^{G\times V}N_e^{G\times V}C
 \simeq (N_e^VC)^{\otimes G}\longrightarrow N_e^VC.
\]
It is projection to the identity $G$-coordinate: this is the
Cartesian counit for the pointed inclusion $BV\to BG\times BV$.
Geometric fixed points preserve finite products of coalgebras and
therefore take it to the same projection for $\Phi^VN_e^VC$.
Using
$x^*\Phi^V\simeq\Phi^V(x\times\id_{BV})^*$ the adjunction
\eqref{eq:point-adjunction} identifies the resulting edge as
Cartesian. Thus $V\star s\in\NormSys_p^{\mathrm{bb}}$.

For the formula on arrows, put $K=\ker(W\to U)$. Norm--inflation
base change identifies the source of the section map as
\[
 N_W^V\infl_U^W N_e^UC\simeq N_K^V\infl_e^K C.
\]
Under the adjunction of Lemma~\ref{lem:lift}, both the section map
and $N_K^V\bar{\Delta}_K$ have adjunct
$C^{\otimes(V/K)}\to C$, given by projection to the factor
indexed by the zero coset $K$. Hence the two maps agree.
Applying $\Phi^V$ gives the induced map $U\star C\to V\star C$.
Under the canonical identifications of both objects with $C$,
this is $\Phi^K\bar{\Delta}_K=\varphi_K$.
\end{proof}

Let $\mathcal W'$ be the backward morphisms in $Q\Vect$.
We write $\Kpart(\Fp)$ for Yuan's partial algebraic
$K$-theory, which is an $\mathbb{E}_\infty$-monoid characterized by the symmetric monoidal equivalence
\[
 Q\Vect[(\mathcal W')^{-1}]\simeq B\Kpart(\Fp)
\]
\cite[Theorem~4.14 and Corollary~4.23.1]{Yuan}.
Here $BM$ denotes the one-object category with endomorphism
anima $M$.
The group completion of $\Kpart(\Fp)$ is ordinary algebraic $K$-theory
$K(\Fp)$ \cite[Corollary~4.9.2]{Yuan}.
Dimension identifies $\pi_0\Kpart(\Fp)$ with $\mathbb N$ and hence
gives a canonical map $\Kpart(\Fp)\to\mathbb N$
\cite[Proposition~4.10]{Yuan}.
We will use the following forward-localization variant.

\begin{lemma}\label{lem:forward-localization}
Let $\mathcal W$ be the class of forward morphisms in $Q\Vect$, induced by inclusions. There is a
monoidal equivalence
\[
 Q\Vect[\mathcal W^{-1}]\simeq B\Kpart(\Fp).
\]
The dimension map $\Kpart(\Fp)\to\mathbb N$ is an
$\Fp$-homology equivalence.
\end{lemma}
\begin{proof}
We combine
duality of finite-dimensional vector spaces with
the canonical isomorphism $Q(\mathcal M^{\op})\cong Q(\mathcal M)$, see e.g. Quillen
\cite[\S2, equation~(5)]{Quillen}. We get 
a covariant monoidal self-equivalence of $Q\Vect$
interchanging $\mathcal W'$ and $\mathcal W$.
Explicitly it takes
$U\twoheadleftarrow W\hookrightarrow V$ to
$U^\vee\twoheadleftarrow U^\vee\times_{W^\vee}V^\vee\hookrightarrow V^\vee$.
It therefore identifies the forward localization with
$Q\Vect[(\mathcal W')^{-1}] \simeq B\Kpart(\Fp)$ as described above.
The last claim about the dimension map is \cite[Theorem~5.2]{Yuan}.
\end{proof}

\begin{proof}[Proof of Theorem~\ref{thm:coalgebra-frobenius-action}]
Proposition~\ref{prop:Q} and Lemma~\ref{lem:forward-localization}
give a monoidal action of $B\Kpart(\Fp)$ on $\coAlg_p^{\mathrm{bb}}$.
Equivalently, this is a map of $\mathbb{E}_2$-monoids from
$\Kpart(\Fp)$ to $\End(\id_{\coAlg_p^{\mathrm{bb}}})$.
The underlying anima of the latter is $\Fp$-homology-local because
it is a limit of coalgebra mapping anima, which are limits of
mapping anima between bounded-below $p$-complete spectra.
The map therefore factors uniquely through the $\Fp$-homology
equivalence $\Kpart(\Fp)\to\mathbb N$.
Since mod-$p$ homology equivalences are preserved by finite products,
this factorization takes place in $\mathbb{E}_2$-monoids.
Thus the action descends uniquely to a $B\mathbb N$-action
(see also \cite[\S6.2]{Yuan}). 
To identify the generator, consider the spans
\[
 a_V=(0\leftarrow0\hookrightarrow V),\qquad
 b_V=(0\leftarrow V\xrightarrow{=}V).
\]
In the forward localization, $a_V$ is inverted and the loop
$a_V^{-1}b_V$ acts by $\varphi_V$. Under the duality used in the proof of
Lemma~\ref{lem:forward-localization}, this is the usual partial
$K$-theory class of $V^\vee$.
Taking $V=\Fp$ identifies the generator $1\in\mathbb N$ with
$\varphi_{\Fp}=\varphi_p$ by \eqref{eq:frob} and
\eqref{eq:coalgebra-frobenius}. This proves the theorem.
\end{proof}


\end{document}